\documentclass{article}
\usepackage{latexsym}
\usepackage[all]{xy}
\usepackage{amssymb,amsmath,amsthm,verbatim,color}
\usepackage{enumitem}

\newtheorem{theorem}{Theorem}
\newtheorem{lemma}[theorem]{Lemma}

\newtheorem{corollary}[theorem]{Corollary}

\title{On the Riemann-Roch formula: old and new}
\author{Claudio Fontanari}
\date{}

\begin{document}

\maketitle

\begin{abstract}
The Riemann-Roch formula is a cornerstone in the classical theory of algebraic curves. Here we 
present a novel approach to its proof, by answering
a question posed in 2007 by Matthew Baker and Serguei Norine.
\vspace{0.01cm}

\vspace{0.1cm}

\end{abstract}

\vspace{0.5cm}

Algebraic curves defined over the complex numbers are a basic but beautiful topic in classical algebraic geometry. Considered as a Riemann surface, an algebraic curve has a natural notion of \emph{genus}, counting the number of holes in the corresponding topological surface. On the other hand, when embedded into a projective space, an algebraic curve is cut out by any linear family of planes into a linear system of \emph{divisors}, each one made by finitely many points lying on the curve (counted with multiplicity). The number of such points is called the \emph{degree} of the divisor and turns out to be constant in the associated linear system. Two divisors in the same linear system are said to be 
\emph{linearly equivalent}. 
A precise relation among the dimension of a linear system, its degree and the genus of the underlying algebraic curve is provided by the celebrated Riemann-Roch Theorem. In order to state it properly, we need to introduce a little piece of notation. 

Let $C$ be a smooth algebraic curve of genus $g$. For every divisor $D$ on $C$, denote by $\deg(D)$ its degree and by $r(D)$ the (projective) dimension of the associated linear system. The Riemann-Roch formula is the following:
$$
r(D) - r(K-D) = \deg(D) - g + 1
$$
where $K$ is a suitable divisor on $C$ of degree $2g-2$, the so-called \emph{canonical divisor}, and addition between divisors is defined in the natural way by formal sums of points. The proof of the Riemann-Roch formula is of course well established, but still is not completely trivial: the standard classical approach is presented in modern terms for instance in the textbook \cite{Ful}, \S 8.6. 

Quite recently, Matthew Baker and Serguei Norine introduced in \cite{BN} an original approach to the Riemann-Roch formula, aimed at its translation from the geometric realm of algebraic curves to the combinatorial world of finite graphs. 

Now we rephase Baker and Norine's statement in the context of classical algebraic geometry. According to 
\cite{BN}, Theorem 2.2 and Remark 2.4(i), 
if $K$ is the canonical divisor of $C$ then the following holds:

\begin{theorem} \emph{(Baker-Norine)} The Riemann-Roch formula:
$$
r(D) - r(K-D) = \deg(D) - g + 1
$$
for every divisor $D$ on $C$ is equivalent to the following pair of properties:

(RR1) If $r(D) < 0$ then there exists a divisor $N$ of degree $g-1$ with $r(N)<0$
such that $r(N-D)\ge 0$.

(RR2) We have $r(K) \ge g-1$.

\end{theorem}

In \cite{BN}, Remark 2.5 (i), the authors point out that \emph{When $C$ is a Riemann surface} (...) \emph{one can show independently of the Riemann-Roch theorem that $r(K)=g-1$} (...) \emph{Thus one can prove directly that (RR2) holds. We do not know if there is a direct proof of (RR1) which does not make use of Riemann-Roch, but if so, one could deduce the classical Riemann-Roch theorem from it using Theorem 2.2}. 

Here we present such a direct proof, thus answering Baker and Norine's interesting question in the affirmative. 

Our tools are as elementary as the proof of (RR2) 
and for a careful exposition we can refer to the textbook \cite{Ful}, \S 8.3 (Riemann's theorem), \S 8.5 (indeed, the last Corollary of this section is precisely (RR2)), and \S 8.6 (Noether's reduction lemma): 

\begin{theorem} \emph{(Riemann's inequality)}
For every divisor $D$ on $C$ we have $ r(D) \ge \deg(D)-g$.
\end{theorem} 

\begin{lemma} \emph{(Noether's reduction)} 
Let $P \in C$ and let $D$ be an effective divisor on $C$. If $r(K-D-P) < r(K-D)$ then $r(D+P) = r(D)$.
\end{lemma}

We point out the following easy consequence of Noether's reduction:

\begin{corollary}\label{NR}
Let $P \in C$ be a general point and let $D$ be a non effective divisor on $C$ of degree $d$. If $d \le g-2$ then $r(D+P) = r(D) = -1$.
\end{corollary}

Indeed, assume by contradiction that $r(D+P) = 0$. 
Since $\deg(D+P) = d+1 \le g-1$, we have
$r(K-D-P) \ge 0$. Hence $r(K-D-P-Q) < r(K-D-P)$ for every general point $Q \in C$ and from Noether's reduction we may deduce 
that $r(D+P+Q) = r(D+P) = 0$, hence $D+P+Q$ is linearly equivalent to a unique effective divisor $E$ passing through $Q$. By exchanging the role of the two general points $P$ and $Q$, we get that $D+Q+P = D+P+Q$ is linearly equivalent to a unique effective divisor, necessarily the same $E$, passing through $P$. Therefore $D+P+Q$ is linearly equivalent to $E = P+Q+R$ with $R$ effective, hence it follows that $D$ is effective, a contradiction.

In order to prove (RR1), let $D$ be a non effective divisor on $C$ of degree $d$. By Riemann's inequality we may assume $d \le g-1$. If $d = g-1$ then (RR1) is obviously satisfied with $N = D$, hence we may assume 
$d \le g-2$. Let now $n = (g-1) - d \ge 1$ and pick $n$ general points $P_1, \ldots P_n$ on $C$. By taking $N = D + P_1 \ldots + P_n$, we obtain that 
$N-D = P_1 + \ldots + P_n$ is effective, $\deg(N) = d + n = g-1$ and $r(N) = r(D + P_1 \ldots + P_n) = 
r(D + P_1 \ldots + P_{n-1}) = \ldots = r(D + P_1) = r(D) = -1$ by applying $n$ times Corollary \ref{NR}.

\vspace{0.5cm}
\noindent \textbf{Acknowledgements.} The author 
discussed Baker and Norine's question with Edoardo Sernesi, who insightfully pointed him towards Noether's reduction Lemma. The author is a member of 
GNSAGA of INdAM. 

\vspace{0.5cm}
\noindent \textbf{Disclosure statement.}
The author confirms that there are no relevant financial or non-financial competing interests to report.

\vspace{0.5cm}

\noindent
Claudio Fontanari \\
Dipartimento di Matematica \\
Universit\`a degli Studi di Trento \\
Via Sommarive 14, 38123 Trento (Italy) \\
claudio.fontanari@unitn.it

\end{document}